# Efficient arc-flow formulations for makespan minimisation on parallel machines with a common server

Alessandro Druetto[1], Andrea Grosso[1], Jully Jeunet[2]*, Fabio Salassa[3]

[1]Dipartimento di Informatica, Università di Torino, Corso Svizzera 185, Torino 10149, Italy

[2]CNRS, Université Paris Dauphine, PSL Research University, CNRS UMR[7243], Lamsade, place du Maréchal de Lattre de Tassigny, 75775 Paris Cedex 16, France

[3]Dipartimento di Ingegneria Gestionale e della Produzione (DIGEP), Politecnico di Torino, Corso Duca degli Abruzzi 24, Torino 10129, Italy

## Abstract

We consider the problem of scheduling non preemptively a set of jobs on parallel identical machines with prior setup operations on a single shared server, where the objective is to minimise the makespan. We develop an arc-flow formulation to the problem with two multigraphs, one for the machines and one for the server, with a same set of nodes representing points in time, and arcs associated with job execution, and with machines or server idleness. The resulting Flow-Flow formulation (FFF) and its tuned version (FFT) are compared with the best existing model in the literature, TIV I, on benchmark instances with up to 200 jobs and 7 machines. Computational results showed that our Flow-Flow models outperformed TIV I for instances with 100 jobs and produced optimal solutions to the vast majority of problems with 150 and 200 jobs for which TIV I failed to even find integer solutions within the standard limited runtime of 3600 s. When non optimal, solutions provided by our models exhibited very low gaps to best bound.

## 1 Introduction

In many manufacturing systems, a common server such as a human operator, a robot or a tool needs to be shared by a number of parallel machines to implement setups or loads. Scheduling problems with a single server occur frequently in automated material handling systems; in flexible manufacturing systems where an automated guided vehicle is used to load jobs on machines (Hall et al., 2000); in the printing industry where a team of workers must clean and reset presses each time a new order is received (Huang et al., 2010) or similarly in knitted fabrics with knitting machines needing to be emptied and their needles, repositioned (Kerkhove and Vanhoucke, 2014). Sharing the server resource results in machine idle time that can be reduced or eliminated by developing a good schedule.

The focus of this paper is on the identical parallel machine scheduling problem with sequence independent setup times and a single shared server where the objective is to minimise the makespan. The problem denoted as $Pm|S1|C_{\max}$ is to schedule a set of jobs $N = \{1, 2, \ldots, n\}$ on an arbitrary number $m$ of identical parallel machines ($m \geq 2$), where each job $j \in N$ must be processed non preemptively on one of the machines for $p_j$ units of time. Prior to processing, $s_j$ units of time must be spent for setup operations on a single shared server. Since the server is shared among all machines, no more than one setup operation can take place at each point in time.

---

*Corresponding author. E-mail addresses: alessandro.druetto@unito.it, andrea.grosso@unito.it, jully.jeunet@dauphine.fr, fabio.salassa@polito.it.

Seminal contributions to the problem mostly consider two machines with equal processing or setup times, and provide complexity analyses and solution methods (Koulamas, 1996; Kravchenko and Werner, 1997, 2001; Glass et al., 2000; Hall et al., 2000 among others). For an extended and clear survey that uses the standard three-field notation described by Graham et al. (1979), the reader may refer to Bektur and Sarac (2019).

The problem with an arbitrary number of machines $m \geq 2$ and general job sets ($\forall s_j, p_j$) was first optimally addressed by Kim and Lee (2012) who provided two Mixed Integer Programming (MIP) formulations. The first one uses the sequence of setups on the server whereas the second one relies on the server waiting time. Experiments from Koulamas (1996) were adapted to create server waiting time for $m \geq 2$ machines since techniques in Abdekhodaee and Wirth (2002) can provide optimal solutions when there is no server waiting time. Instances with $n = 10$ jobs were exactly solved in less than one minute. For larger problems up to 40 jobs and 6 machines, 40 to 70% of optimal solutions on average were found within the limited runtime of 3600 s. Elidrissi et al. (2018a) addressed the same scheduling problem for which they proposed two MIP formulations based on completion time variables and time-indexed variables with better performance. The authors then enhanced these formulations with strengthening constraints (Elidrissi et al., 2021). In addition, they developed three other MIP formulations using respectively network variables, linear ordering variables and completion time variables. The performance of these models was compared with the two MIP formulations of Kim and Lee (2012) on a set of instances generated in a similar way. Results showed that only the time-indexed variable formulation (TIV I) was able to find optimal solutions to some of the instances with more than 10 jobs and up to 100 jobs.

Metaheuristic solution approaches to the problem are available in Kim and Lee (2012) who developed a Simulated Annealing (SA) algorithm combined with tabu search, and in Elidrissi et al. (2020) who proposed a Variable Neighbourhood Search (VNS) algorithm. Elidrissi et al. (2018b) generalised the heuristics of Abdekhodaee and Wirth (2002) and that of Hasani et al. (2016) to the case of an arbitrary number of machines, thus providing two greedy heuristics aimed at minimising machine idle time and server waiting time.

Besides, several variants of the parallel machine scheduling problem with common servers have been considered in the literature, depending upon the objective to minimise, the type of setup times, the assumption about preemption or the number of servers. Changing only one assumption of the problem considered in this paper, we get the following contributions.

Liu et al. (2019) developed a branch-and-bound algorithm to minimise the weighted job completion time which was able to optimally solve instances with no more than 20 jobs and 3 machines. Abu-Shams et al. (2022) proposed a heuristic-based Genetic Algorithm (GA) to minimise tardiness so as to deal with large-sized problems up to 2000 jobs and 10 machines.

Sequence-dependent setup times are considered in Hamzadayi and Yildiz (2017) who presented a mixed integer linear programming (MILP) model for small-sized instances (no more than 20 jobs and 5 machines) as well as a SA and a GA for larger problems (up to 100 jobs and 10 machines). Silva et al. (2021) proposed an arc-time-indexed formulation able to solve exactly larger instances (21 jobs, 7 machines).

Cheng et al. (2017) presented complexity analyses to the problem with preemption whereas Elidrissi et al. (2022) provided a MIP model to deal with 2 servers and small problems (10 jobs, 5 machines) as well as a VNS for large instances (250 jobs, 5 machines), but for regular jobs only (regular jobs are such that $p_i \leq p_j + s_j, \forall (i,j)$).

Existing solution approaches to the scheduling problem with common servers are thus barely able to provide optimal solutions to instances with more than 20 jobs in a fast compu-tation time. We propose an arc-flow formulation for problem $Pm|S1|C_{\max}$ which is capable to outperform the best model of Elidrissi et al. (2021), namely their improved time-indexed variables formulation (TIV I), both in terms of number of optimal solutions and computation time. While their formulation could optimally solve half of a quite limited number of instances with 100 jobs within 3600 s, ours managed to find optimal solutions to an extended number of problems with 100 jobs in a fast computation time. In addition, our models were able to exactly solve the vast majority of instances with 150 and 200 jobs and produced very low gaps to best bound when non optimal.

The remainder of this paper is organised as follows. Our arc-flow formulation is presented in Section 2. Results of the computational experiments are discussed in Section 3. Finally, concluding remarks and directions for future research are provided in Section 4.

## 2 Arc-flow formulation and existing mathematical model TIV I

The proposed arc-flow formulation is presented and illustrated in Section 2.1. Next, we provide a comparison with the model of Elidrissi et al. (2021), TIV I, and we show that the continuous relaxation of our formulation provides a better lower bound than a trivial bound (Section 2.2). Finally we give in Section 2.3 the bounds on the makespan that we used in our arc-flow formulation for which a tuned version is then proposed.

Arc-flow formulations allow for the use of a pseudo-polynomial number of variables and constraints and have been recently applied to classical optimisation problems such as the cutting-stock problem (Martinovic et al., 2018), the bin-packing problem (Brandao et al., 2016) or the berth allocation problem (Kramer et al., 2019a). In the area of scheduling and most closely related to our problem, Mrad and Souayah (2018) proposed an arc-flow formulation for makespan minimisation on identical parallel machines and showed its efficiency to solve most of the hard instances from the literature. Gharbi and Bamatraf (2022) provided an improved arc-flow model for the same problem, with enhanced bounds. Results on benchmark instances with up to 200 jobs and 100 machines showed the superiority of their model over that of Mrad and Souayah (2018). Kramer et al. (2019b) also considered the scheduling problem on identical parallel machines but with the aim of minimising the total weighted completion time. They developed enhanced arc-flow formulations able to solve exactly large-sized instances up to 400 jobs. This work was then extended in Kramer et al. (2020) to jobs with release dates. To the best of our knowledge, arc-flow formulations have not been used to address the scheduling problem on identical parallel machines with a common server.

### 2.1 The Flow-Flow Formulation (FFF)

We first provide the multigraph representation of the problem and we introduce the notations. The problem modelling is then illustrated with a scheduling example. Finally, we present the mathematical formulation and the procedure to handle identical jobs.

#### 2.1.1 Multigraph representation and notations

Our arc-flow formulation uses two multigraphs $G_{K=\{M,S\}}(V, A_K)$ in order to model the scheduling of the $m$-machines collection, $M$, and the single server $S$, respectively. We therefore call our model Flow-Flow formulation (FFF in the following), having an arc-flow formulation for both machines and server. The two graphs have the same set of nodes $V$, each node representing a unit time slot

$$V = \{0, 1, 2, \ldots, T\},$$

where $T$ is the time horizon over which the set of jobs $N = \{1, 2, \ldots, n\}$ must be scheduled. The set $A_K$ of arcs of multigraph $K = \{M, S\}$ is defined as

$$
\begin{aligned}
A_K &= \{a_{jt}=(t,t+b_j;j)\colon j\in N, t=0..T-s_j-p_j\} \cup \{a_t=(t,t+1)\colon t=0..T-1\}, \\
\text{with } b_j &= \begin{cases} s_j+p_j, & K=M, \\ s_j, & K=S. \end{cases}
\end{aligned}
$$

An arc $a_{jt}$ represents the possible execution of job $j$ in time interval $(t, t+b_j)$ so it links node $t$ (start time) to node $t+b_j$ (end time) both on the machines graph $G_M$ and on the server graph $G_S$. An arc $a_t$ expresses the idleness of the server or machines in time interval $(t, t+1)$. In addition, we let $\delta^-_{[G_K,t]}$, $\delta^+_{[G_K,t]}$ be respectively the set of ingoing and outgoing arcs $a_{jt}$ in/from node $t$ in graph $G_K$, $K=\{M,S\}$.

We let $x_{jt}$ be a binary variable that takes a value of 1 if job $j$ starts at time $t$ and 0 otherwise. If $x_{jt}=1$, this unit flow is placed on arc $a_{jt}$ that links node $t$ to node $t+s_j$ in graph $G_S$ and node $t$ to node $t+s_j+p_j$ in graph $G_M$. Job $j$ thus starts at time $t$ both on the server and on a machine since the non preemption assumption leads to reserve a machine for processing the job as soon as its setup starts on the server. Thus, the machine is not considered as idle on $t$ if the setup on the server starts on $t$.

We let $y_t^M$, $y_t^S$ be the integer variables expressing the number of idle machines and server respectively. Variable $y_t^S$ is obviously binary since we consider a single server. The flow value $y_t^M$ (resp. $y_t^S$) is placed on arc $a_t$ in graph $G_M$ (resp. $G_S$) that connects node $t$ to node $t+1$.

Finally we introduce a binary variable $z_t$ that takes a value of 1 if the last scheduled job ends at time $t$ and 0 otherwise. If $z_t=1$ the makespan is therefore equal to $t$.

Table 1 summarises the notations and definitions.

| | |
|---|---|
| **Indices** | |
| $j$ | Job |
| $t, \tau$ | Time |
| $K=\{M,S\}$ | Type of multigraph ($M$ for the machines, $S$ for the server) |
| **Parameters** | |
| $T$ | Time horizon (set to some upper bound on the makespan) |
| $N=\{1,2,\ldots,n\}$ | Set of jobs to be scheduled |
| $s_j$ | Setup time of job $j$ (on the server) |
| $p_j$ | processing time of job $j$ (on a machine) |
| **Definitions** | |
| $G_{K=\{M,S\}}(V, A_K)$ | Multigraph for $K=\{M,S\}$ |
| $V=\{0,1,2,\ldots,T\}$ | Set of nodes shared by both multigraphs |
| $a_{jt}=(j;t,t+b_j)$ | Execution arc of job $j$ from $t$ to $t+b_j$; $b_j=s_j$ in $G_S$, $b_j=s_j+p_j$ in $G_M$, $j\in N$, $t=0..T-s_j-p_j$. |
| $a_t=(t,t+1)$ | Idleness arc from $t$ to $t+1$, $t=0..T-1$ |
| $A_K=\{a_{jt}\}\cup\{a_t\}$ | set of arcs in the multigraph $K=\{M,S\}$ |
| $\delta^-_{[G_K,t]}$ | set of ingoing execution arcs $a_{jt}$ in node $t$ in $G_K$, $K=\{M,S\}$ |
| $\delta^+_{[G_K,t]}$ | set of outgoing execution arcs $a_{jt}$ from node $t$ in $G_K$, $K=\{M,S\}$ |
| **Variables** | |
| $x_{jt}$ | Binary equal to 1 if job $j$ starts at time $t$ and 0 otherwise, $\forall j\in N, \forall t=0..T-s_j-p_j$. |
| $y_t^K$ | Number of idle machines ($K=M$) or server ($K=S$) at time $t$, $\forall t=0..T$ |
| $z_t$ | Binary equal to 1 if the last scheduled job ends at time $t$, $\forall t=1..T$ |

Table 1: Notations and definitions

### 2.1.2 Illustration of the problem modelling

To illustrate the problem modelling, let us consider $n = 5$ jobs with $\{s_1, \ldots, s_5\} = \{2, 3, 3, 2, 2\}$ and $\{p_1, \ldots, p_5\} = \{3, 5, 4, 5, 3\}$. The jobs are to be scheduled on a time horizon $T = 18$ with one server and $m = 3$ machines. Figure 1 displays on a same time scale the Gantt chart (upper part) of a feasible solution to the problem with $C_{\max} = 17$ and the associated flows on multigraphs $G_M$ and $G_S$ (lower part). The flows are generated at node 0 for both graphs. For instance, job 1 starts at time 0 ($x_{1,0} = 1$) so an execution arc with a unit flow is placed on $G_S$ from node 0 to node $0 + s_1 = 2$ and the corresponding nodes $\{0, 1, 2\}$ are connected with zero-flow (hence not drawn) idleness arcs since the single server is busy. On graph $G_M$ an execution arc connects node 0 to node $0 + s_1 + p_1 = 5$ even if the processing of the jobs starts on $M_1$ only at time 2, as the machine is reserved as soon as the setup starts on the server. Thus, between time 0 and 5, $M_1$ is not idle; only machines $M_2$ and $M_3$ are idle between time 0 and 2, so idleness arcs with a flow of 2 units connect nodes $\{0, 1, 2\}$. At $t = 2$ job 2 starts its setup on the server, so $M_2$ is reserved and the idleness arc from node 2 to node 3 now carries a unit flow as only $M_3$ is idle. Let us note that in time interval $(8, 10)$ all the machines are busy, so the corresponding nodes are connected with zero-flow idleness arcs. All in all, 3 units of flow are routed on $G_M$, and 1 unit on $G_S$. The sink node is 17 as we have $z_{17} = 1$ so arcs are not drawn between nodes 17 and 18. The flow does not give explicitly an assignment jobs-machines, but it can be decomposed into 3 paths that correspond to machine schedules.

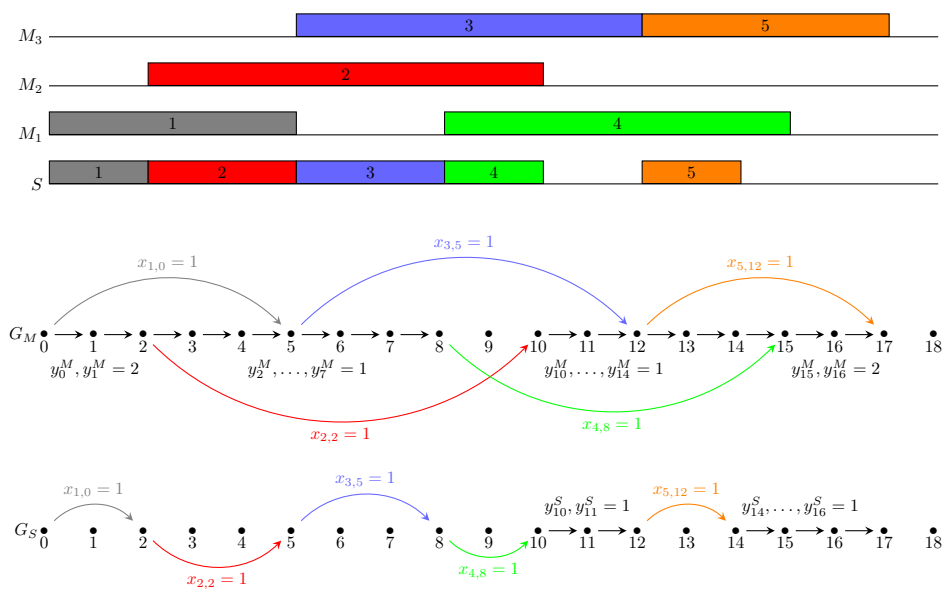

Figure 1: A scheduling example and corresponding flows on $G_M$ and $G_S$.

### 2.1.3 Mathematical formulation (FFF) and identical jobs handling

In order to facilitate the reading of our formulation, we index all equations with time $t$ and we indicate in parenthesis the values of $t$ for which the corresponding constraints hold. Our arc-flow model (FFF) is written as follows

$$\min \sum_{t=1}^{T} t z_t \tag{1}$$

s.t.

$$\sum_{t=0}^{T-s_j-p_j} x_{jt} = 1 \quad (j \in N) \tag{2}$$

$$\sum_{a_{jt} \in \delta^+_{[G,t]}} x_{jt} + y_t^M = m \quad (t=0) \tag{3}$$

$$\sum_{a_{jt} \in \delta^+_{[G,t]}} x_{jt} - \sum_{a_{j\tau} \in \delta^-_{[G,t]}} x_{j\tau} + y_t^M - y_{t-1}^M = -m z_t \quad (t=1..T-1) \tag{4}$$

$$- \sum_{a_{j\tau} \in \delta^-_{[G,t]}} x_{j\tau} - y_{t-1}^M = -m z_t \quad (t=T) \tag{5}$$

$$\sum_{a_{jt} \in \delta^+_{[G_s,t]}} x_{jt} + y_t^S = 1 \quad (t=0) \tag{6}$$

$$\sum_{a_{jt} \in \delta^+_{[G_s,t]}} x_{jt} - \sum_{a_{j\tau} \in \delta^-_{[G_s,t]}} x_{j\tau} + y_t^S - y_{t-1}^S = -z_t \quad (t=1..T-1) \tag{7}$$

$$- \sum_{a_{j\tau} \in \delta^-_{[G_s,t]}} x_{j\tau} - y_{t-1}^S = -z_t \quad (t=T) \tag{8}$$

$$x_{jt} \in \{0,1\} \quad (j \in N, t=0..T-s_j-p_j) \tag{9}$$

$$y_t^M \geq 0, \quad y_t^S \in \{0,1\} \quad (t=0..T-1) \tag{10}$$

$$z_t \in \{0,1\} \quad (t=1..T) \tag{11}$$

Constraints (2) require that a starting time is assigned to each job $j \in N$. Constraints (3)–(5) related to the machines make variables $x_{jt}, y_t^M$ a set of flow variables, requiring $m$ units of flow to be routed on graph $G_M$ from source node 0 to some sink node $t \in \{1, \ldots, T\}$ for which $z_t = 1$. Let us note that summing up constraints (3)–(5) we get $m(1 - \sum_{t=1}^{T} z_t) = 0$, so at most one $z_t$ can take a value of 1.

Constraints (6)–(8) make variables $x_{jt}, y_t^S$ a set of flow variables, requiring one unit of flow to be routed on graph $G_S$ from node 0 to some node $t \in \{1, \ldots, T\}$ for which $z_t = 1$. These constraints guarantee that no more than one job can be processed on the server simultaneously. Given that $\sum_{t=1}^{T} z_t = 1$ must hold in all feasible solutions, the objective function (1) correctly represents the makespan.

In sets of jobs to be scheduled, it is not uncommon that quite large families of identical jobs — with identical setup and processing times — emerge. In such cases, a powerful optimisation procedure consists in keeping a single copy of identical jobs and schedule it a number of times equal to the number of copies. Formally, letting $n_j$ be the number of jobs identical to $j$ in set $N$, we keep only one copy of such a job $j$, and schedule it $n_j$ times. We thus eliminate all duplicated jobs from $N$, so that no two jobs in $N$ are identical; we compute the multiplicity $n_j$ of each job $j$ and replace constraints (2) with

$$\sum_{t=0}^{T-s_j-p_j} x_{jt} = n_j.$$

This procedure presents similarities with the identical job grouping adopted in Kramer et al. (2019b). In many instances the procedure can substantially reduce the number of variables involved in the model. Due to its simplicity we always apply it to the FFF model.

## 2.2 Comparison with the timed-indexed formulation of Elidrissi et al. (TIV I)

The time-indexed formulation of Elidrissi et al. (2021), TIV I, makes use of the same binary variables $x_{jt}$, and an auxiliary variable $C_{\max}$ so as to implement a classical min-max formulation

$$\min C_{\max} \tag{12}$$

s.t.

$$\sum_{t=0}^{T} x_{jt} = 1 \qquad \forall j \in N \tag{13}$$

$$\sum_{j=1}^{n} \sum_{s=(t-s_j-p_j+1)^+}^{t} x_{jt} \leq m \qquad \forall t = 0..T \tag{14}$$

$$\sum_{j=1}^{n} \sum_{s=(t-s_j+1)^+}^{t} x_{jt} \leq 1 \qquad \forall t = 0..T \tag{15}$$

$$\sum_{t=0}^{T} (t + s_j + p_j) x_{jt} \leq C_{\max} \qquad \forall j \in N \tag{16}$$

$$x_{jt} \in \{0,1\} \qquad \forall j \in N \quad \forall t = 0..T \tag{17}$$

Constraints (13) require that each job is assigned a starting time. Constraints (14) ensure that at each point in time no more than $m$ jobs are executed on the parallel machines. Constraints (15) guarantee that one job at most is processed by the server at time $t$. Constraints (16) force $C_{\max}$ to equal the completion time of the last scheduled job.

TIV I model expressed by Eq. (12)–(17) has quite a weak relaxation, often providing a lower bound below a trivial one, $\text{LB}_{\text{Trivial}}$, given by

$$\text{LB}_{\text{Trivial}} = \frac{1}{m} \sum_{j=1}^{n} (s_j + p_j). \tag{18}$$

Empirically, for all the 300 instances on which TIV I was applied in in our numerical experiments it actually turned out that $\text{LB}_{\text{Trivial}}$ was much greater than the lower bound from the continuous relaxation of TIV I (see Section 3). On the other hand, our Flow-Flow model FFF in Eq. (1)–(11) exhibits the following property.

**Property 1.** *The optimal value of the continuous relaxation of the Flow-Flow model (FFF) is never less than $\text{LB}_{\text{Trivial}}$.*

*Proof.* We consider any feasible solution of the continuous relaxation and we prove that its objective function value cannot be less than the trivial lower bound. In order to ease the notation, we let $I_t = \sum_{a_{j\tau} \in \delta^-_{[G_M,t]}} x_{jt}$ and $O_t = \sum_{a_{jt} \in \delta^+_{[G_M,t]}} x_{jt}$, with $I_0 = 0$ and $O_T = 0$, so that constraints (3)–(5) can be rewritten, changing sign, as

$$-m = I_t - O_t \quad\quad\quad - y_t^M \quad\quad (t=0)$$
$$m z_t = I_t - O_t \quad + y_{t-1}^M - y_t^M \quad\quad (t=1..T-1)$$
$$m z_t = I_t - O_t \quad + y_{t-1}^M \quad\quad\quad (t=T).$$

We multiply each constraint by the corresponding $t$, sum them up and divide by $m$

$$\sum_{t=1}^{T} t z_t = \frac{1}{m} \sum_{t=0}^{T} t[I_t - O_t] + \frac{1}{m}\left[\sum_{t=1}^{T-1}(t y_{t-1}^M - t y_t^M) + T y_{T-1}^M\right].$$

We note the following.

i $\sum_{t=1}^{T} t z_t$ is precisely the objective function.

ii The sum $\sum_{t=0}^{T} t[I_t - O_t]$ is a weighted sum of all arc-flow variables $x_{jt}$ of arcs $a_{jt}$. Hence, every $x_{jt}$ will appear in the sum exactly twice: the first time in $O_t$, for the tail of arc $a_{jt}$, multiplied by $-t$, and the second time in $I_{t+s_j+p_j}$, for the head of the arc, multiplied by $(t + s_j + p_j)$. Hence the contribution of $x_{jt}$ to the sum is $(s_j + p_j)x_{jt}$. Thus

$$\sum_{t=0}^{T} t[I_t - O_t] = \sum_{j=1}^{n} \sum_{t=0}^{T-s_j-p_j} (s_j + p_j) x_{jt}$$

and by constraint (2)

$$\sum_{t=0}^{T} t[I_t - O_t] = \sum_{j=1}^{n} \left[ (s_j + p_j) \sum_{t=0}^{T-s_j-p_j} x_{jt} \right] = \sum_{j=1}^{n} (s_j + p_j).$$

iii For the sum involving the $y_t$ variables,

$$\sum_{t=1}^{T-1} (t y_{t-1}^M - t y_t^M) = \sum_{t=1}^{T-1} t y_{t-1}^M - \sum_{t=1}^{T-1} t y_t^M =$$

$$= \sum_{t=0}^{T-2} (t+1) y_t^M - \sum_{t=1}^{T-1} t y_t^M =$$

$$= y_0^M + \sum_{t=1}^{T-2} y_t^M - (T-1) y_{T-1}^M$$

hence

$$\frac{1}{m} \left[ \sum_{t=1}^{T-1} (t y_{t-1}^M - t y_t^M) + T y_{T-1}^M \right] = \frac{1}{m} \sum_{t=0}^{T-1} y_t^M$$

By the above considerations, we get

$$\sum_{t=1}^{T} t z_t = \frac{1}{m} \sum_{j=1}^{n} (s_j + p_j) + \frac{1}{m} \sum_{t=0}^{T-1} y_t^M,$$

showing that the objective function is made of the trivial lower bound $\mathsf{LB}_{\text{Trivial}}$ in Eq. (18) plus a nonnegative contribution from the idle variables $y_t^M$, thus proving the claim. □

This property establishes a minimum performance for our Flow-Flow Formulation but there is no proof of dominance of the lower bound from its continuous relaxation over that of TIV I.

## 2.3 Bounds on the makespan and tuned version of the Flow-Flow model

We first provide the lower and upper bounds on the makespan that we use in our Flow-Flow models. We then detail the features of the tuned version of the Flow-Flow formulation.

### 2.3.1 Lower and upper bounds on the makespan

As a better lower bound than $\text{LB}_{\text{Trivial}}$ on the optimal makespan $C_{\max}^{\text{opt}}$, we use the same as Elidrissi et al. (2021) who state that (i) if there is no server waiting time in the optimal schedule, then $C_{\max}^{\text{opt}}$ is equal to the sum of setup times and the shortest processing time; (ii) if there is no machine idle time, then $C_{\max}^{\text{opt}}$ is equal to the average over the number of machines of the sum of all jobs execution times, and the weighted sum of the $(m-1)^{\text{th}}$ first jobs setup times ranked in increasing order $\{\sigma(j)\}_{j=1..n}$. Formally this lower bound, $\text{LB}_{\text{Better}}$, is given by

$$\text{LB}_{\text{Better}} = \max\left(\sum_{j=1..n} s_j + \min_{1\leq j\leq n} p_j, \text{LB}_{\text{Trivial}} + \frac{\sum_{j=1}^{m-1}(m-j)s_{\sigma(j)}}{m}\right). \tag{19}$$

To derive an upper bound on the optimal makespan, Elidrissi et al. (2021) used the two greedy heuristics HS1, HS2 they developed in Elidrissi et al. (2018b). Heuristic HS1 aims at minimising machines idle time whereas HS2 seeks to minimise the server idle time. Both heuristics are based on ordering jobs according to six priority rules, namely shortest or longest processing times (SPT, LPT), setup times (SST, LST) and completion times (SCT, LCT). As many jobs in a list can have identical processing or setup times, we introduce a tie breaking rule that consists of arranging jobs with same criterion value according to a second criterion consistent with the first one. For instance, if rule LPT is used, jobs with same processing times in the ranked list are secondarily ordered according to LST. Analogously, with SPT we apply SPT+SST. Finally for SCT, we use SCT+SPT. The upper bound UB provides the value of the time horizon length $T$ in models FFF and TIV I, and is formally defined as

$$T = \text{UB} = \min\left(C_{\max}^{\text{HS1}}, C_{\max}^{\text{HS2}}\right). \tag{20}$$

### 2.3.2 Tuned version of the Flow-Flow model (FFT)

The Flow-Flow Tuned model (FFT) is also defined by Eq. (1)–(11) as FFF, but includes several features that we list below.

**Formulation strengthening**. To strengthen the formulation, we use the lower bound $\text{LB}_{\text{Better}}$ of Elidrissi et al. (2021) in Eq. (19) by setting to zero variables $z_t$ in Eq. (10) for all $t < \text{LB}_{\text{Better}}$.

**Branching priority and direction**. Branching priority is enforced for variables $x_{jt}$ during the branch-and-cut phase, with decreasing priority over the increase of $t$. The branching direction is set to up, so the up branch is taken first at each node, since the aim is to set as soon as possible to 1 the variables close to the start of the sequence. In doing so, we obtain a kind of 'schedule from the beginning' way of branching, which usually makes sense in scheduling, instead of operating on some (almost) randomly selected fractional $x_{jt}$.

**CPLEX parameters configuration**. While investigating the impact of tuned software parameters on the performance, some contributions evidenced that improved parameter configurations may lead to substantial speedup for solving many combinatorial problems (Baz et al., 2007; Hutter et al., 2009 or Fawcett and Hoos, 2016). Pilot runs of our model on test instances led to the configuration of the CPLEX MIP solver for FFT as displayed in Table 2 whereas default parameter values were used for FFF. Table 2 provides the rationale for choosing some specific parameters values. Apart from the branching priority on the variables, the chosen CPLEX parameters values reflect an aggressive setting in searching for good feasible solutions.

| Parameter name | Value | Description/Motivation |
| --- | --- | --- |
| MIP dynamic search switch | 1 | Traditional branch-and-cut, no Dynamic Search for the branching phase. |
| Feasibility pump switch | 2 | Focus on finding solutions with better objective values, instead of potentially worse feasible ones (for further details on the method, see Fischetti et al., 2005). |
| RINS heuristic frequency | 3 | Application every 3 nodes of the Relaxation Induced Neighborhood Search heuristic, an expensive heuristic useful to find high quality integer solutions (this heuristic is presented in Danna et al., 2005). |
| MIP probing level | 2 | Enforcement of an aggressive probing on variables before the branching phase. |
| MIP priority order switch | 1 | Required for the aforementioned branching priority rule, otherwise it would be ignored. |
| MIP repeat presolve switch | 3 | Repetition of the presolve to allow new cuts and new root cuts. |
| MIP dive strategy | 3 | The MIP traversal strategy occasionally performs probing dives, where it looks ahead before deciding which node to choose. With a value of 3, the solver is enabled to spend more time exploring potential solutions that are similar to the current one. |

Table 2: CPLEX configuration for FFT

## 3 Computational experiments

The exact models TIV I, FFF and FFT were coded in C++ and linked with the CPLEX Callable Library (C API) version 22.1. Computational experiments were executed in Ubuntu Linux 22.04 under WSL2, on an Intel(R) Core(TM) i7-10700 CPU with 2.90GHz and 16GB of RAM machine, running Windows 11. We first present the benchmark instances. Results are then discussed, starting with small-sized problems ($n = 10, 20, 50$ jobs), proceeding with medium-sized ones ($n = 100$ jobs) and ending with the largest instances ($n = 150, 200$ jobs).

### 3.1 Benchmark instances

The first set of instances we used are those of Elidrissi et al. (2021) who made them publicly available at https://sites.google.com/site/dataforpssproblem/home. We considered only general job sets since, as pointed by Kim and Lee (2012), regularity of all jobs to be scheduled rarely happens in practice. Following Kim and Lee (2012), Elidrissi et al. (2021) generated their instances using four parameters $(n, m, \alpha, \rho)$, with $n$ between 10 and 100 jobs and $m$ varying from 3 to 5 machines. Parameter $\alpha$ determines the interval from which processing times are uniformly drawn as integer values

$$p_j \sim U\left[(1-\alpha)E(p_j), (1+\alpha)E(p_j)\right], \quad E(p_j) = 25, \quad \alpha = \{0.1, 0.3, 0.5\}. \tag{21}$$

Obviously, the smaller $\alpha$ the higher the number of jobs with identical processing times.

Parameter $\rho$ together with $m$ determines the server load $\rho/m$. In order to create idle times for the server, setup times must be lower than processing times. This requirement is satisfied using $E(s_j) = (\rho/m) \cdot E(p_j)$ with $\rho \leq 1$. Setup times are therefore drawn as follows

$$s_j \sim U\left[(1-\alpha)(\rho/m)E(p_j), (1+\alpha)(\rho/m)E(p_j)\right], \quad \rho = \{0.5, 0.7, 1\}. \tag{22}$$

Again, a low $\alpha$ value leads to a high number of jobs with identical setup times.

Elidrissi et al. (2021) considered a subset of 23 combinations of values for parameters $(n, m, \alpha, \rho)$ with 10 replications of random draws for each combination, resulting in a total of 230 instances, 120 of which with $n = 10$ jobs. There are 10 instances with the highest number of jobs $n = 100$ to be scheduled on $m = 5$ machines, which corresponds to 10 draws of $(p_j, s_j)_{j \in N}$ for a single combination $(\alpha, \rho)$. For $(n, m) = (100, 5)$ we extended the number of combinations $(\alpha, \rho)$ to the ones the authors utilised for $n = 50$ and we also considered $m = \{3, 7\}$. The same setting was used with $n = 150, 200$ as shown in Table 3 that summarises the experimental framework where a cross corresponds to the instances of Elidrissi et al. (2021) and a circle indicates our additional instances. In this way, we generated 230 additional instances with $n \geq 100$, having thus a total number of $230 + 230 = 460$ problems to solve.

| | | $n$ | 10 | 20 | 50 | | 100 | | | 150, 200 | |
|---|---|---|---|---|---|---|---|---|---|---|---|
| $\alpha$ | $\rho$ | $m$ | 3, 4 | 3 | 3 | 7 | 3 | 5 | 7 | 3 | 5 | 7 |
| 0.1 | 0.5 | | × | × | × | | ○ | × | | ○ | ○ | |
| | 0.7 | | | | | × | | ○ | ○ | | ○ | ○ |
| | 1.0 | | × | × | | | | | | | | |
| 0.3 | 0.5 | | × | × | | | | | | | | |
| | 0.7 | | × | × | | | | | | | | |
| | 1.0 | | | | | × | | ○ | ○ | | ○ | ○ |
| 0.5 | 0.5 | | | | | × | | ○ | ○ | | ○ | ○ |
| | 0.7 | | × | × | | | | | | | | |
| | 1.0 | | × | × | | | | | | | | |
| # per $(n,m)$ | | | 60 | 60 | 20 | 20 | 20 | 40 | 20 | 20 | 40 | 20 |

Table 3: Combinations of parameters $(n, m, \alpha, \rho)$ in the experiments

## 3.2 Results

The 120 smallest instances with $n = 10$ jobs were all solved to optimality by the 3 exact approaches FFF, FFT and TIV I in a quite fast computation time, less than 0.23 s on average over all instances (0.20 s for FFT and 0.23 s for both FFF and TIV I). For the instances with a higher number of jobs we considered the following indicators.

- #O   Number of Optimal solutions (over 10 draws).
- #N   Number of instances for which No integer solution was found within a time limit of 3600 s (over 10 draws).
- CPU   Average solution time in s, over the number of optimal solutions out of 10 draws.
- DEV$_{\text{CR}}$   Average deviation of optimal solutions to the lower bounds provided by the Continuous Relaxation of each model in %.
- GAP$_{\text{BB}}$   Average gap to Best Bound as provided by CPLEX over non optimal solutions in %.

Medium-sized instances with $n = 20, 50$ were again all solved to optimality by the 3 models but in a much faster computation time by our Flow-Flow formulations compared with TIV I, as shown in Table 4. Obviously, we reported in this table only indicators CPU and $\text{DEV}_{\text{CR}}$ since all solutions were optimal. FFT exhibits the fastest solution time, with about 5 s on average overall cases, and is closely followed by FFF with 36 s on average, which is about 10 times less than TIV I average execution time (355 s). The benefit of using the tuned version of our Flow-Flow formulation FFT over FFF is especially noticeable for the hardest instance $(n, m, \alpha, \rho) = (50, 3, 0.5, 0.5)$ in blue font with an average solution time of about 30 s for FFT, whereas FFF takes about 4 min, a time that is nonetheless quite reasonable compared with the 37 min taken by TIV I on average to find the 10 optimal solutions. The good performance of the Flow-Flow formulations can be explained by the quality of the lower bounds provided by the continuous relaxation, with an average deviation of optimal solutions to these lower bounds ( $\text{DEV}_{\text{CR}}$ ) of 0.46% in all cases, whereas this average deviation equals 75.3% for TIV I. Let us note that indicator $\text{DEV}_{\text{CR}}$ for FFT is even zero in half of the instances, meaning that the continuous relaxation already provides the optimal solutions.

| | | | | **FFT** | | **FFF** | | **TIV I** | |
|---|---|---|---|---|---|---|---|---|---|
| $n$ | $m$ | $\alpha$ | $\rho$ | CPU | $\text{DEV}_{\text{CR}}$ | CPU | $\text{DEV}_{\text{CR}}$ | CPU | $\text{DEV}_{\text{CR}}$ |
| 20 | 3 | 0.1 | 0.5 | 0.60 | 0.43 | 0.82 | 1.31 | 3.34 | 74.00 |
| | | | 1.0 | 1.12 | 1.21 | 1.71 | 2.05 | 6.32 | 73.38 |
| | | 0.3 | 0.5 | 1.24 | 0.00 | 2.18 | 0.23 | 6.34 | 72.33 |
| | | | 0.7 | 1.55 | 0.00 | 3.05 | 0.22 | 11.18 | 71.14 |
| | | 0.5 | 0.7 | 1.75 | 0.00 | 3.26 | 0.21 | 13.62 | 71.14 |
| | | | 1.0 | 6.77 | 0.05 | 78.77 | 0.27 | 97.97 | 70.17 |
| 50 | 3 | 0.1 | 0.5 | 1.98 | 0.07 | 8.10 | 0.15 | 876.22 | 87.85 |
| | | 0.5 | 0.5 | 31.13 | 0.00 | 222.59 | 0.09 | 2207.00 | 87.36 |
| | 7 | 0.1 | 0.7 | 0.36 | 0.95 | 0.84 | 1.57 | 52.86 | 73.86 |
| | | 0.3 | 1.0 | 5.03 | 0.00 | 40.98 | 0.40 | 275.16 | 71.49 |

Table 4: Performance of FFT, FFF and TIV I for $n = 20, 50$ jobs

Table 5 displays the results for the instances with $n = 100$ jobs, applying again the 3 exact approaches. Only FFT was able to optimally solve all the 80 problems whereas FFF failed in 8 cases but provided an integer solution to 6 of them with an extremely low $\text{GAP}_{\text{BB}}$ of about 0.33% (no integer solution could be found within 3600 s for 2 cases only). By contrast, TIV I was unable to achieve optimality for 74% of the problems, and always found low quality integer solutions with $\text{GAP}_{\text{BB}}$ around $33 - 36\%$ after a running time of 3600 s. When TIV I found optimal solutions (21 cases over 80), the average execution time was 37 min, whereas it took respectively 3 s and 17 s on average for FFT and FFF to optimally solve all instances of the same parameters combinations (in black font in Table 5).

Hard instances (in blue font in Table 5) are those for which TIV I found no optimal solution, whereas FFT reached the optimum for all of them in less than 6 min on average. Amongst these 40 hard problems, FFF found the optimal solutions to 30 of them in 26 min on average. Here again, we evidence the benefit of using a tuned version of our Flow-Flow formulation. The hardest combination $(m, \alpha, \rho) = (3, 0.5, 0.5)$ is the one for which TIV I could not even find integer solutions to the 10 draws within 3600 s (#N = 10). This number #N dropped to 2 only with FFF (and to 0 for FFT).

Let us note again the high quality of the lower bound produced by the continuous relaxation of the Flow-Flow models with $\text{DEV}_{\text{CR}}$ that never exceed 0.76% and equal 0 in 50% of the problems for FFT. By contrast, $\text{DEV}_{\text{CR}}$ for TIV I are quite high, about 90% on average, meaning that optimal makespan values are almost twice the LB values.

Referring to the quality of lower bounds obtained from the continuous relaxation, it is worth mentioning that over all the 300 instances with $n \leq 100$, TIV I produced bounds that were 37% lower than $\text{LB}_{\text{Trivial}}$ on average (with a standard deviation of 9.5%) while on the contrary those of FFF were 2.8% greater than $\text{LB}_{\text{Trivial}}$ (with a standard deviation of 2.4%). This result substantiates the proof provided in Section 2.2.

| | | | **FFT** | | **FFF** | | | | **TIV I** | | | | |
|---|---|---|---|---|---|---|---|---|---|---|---|---|---|
| $m$ | $\alpha$ | $\rho$ | CPU | $\text{DEV}_{\text{CR}}$ | #O | #N | CPU | $\text{DEV}_{\text{CR}}$ | $\text{GAP}_{\text{BB}}$ | #O | #N | CPU | $\text{DEV}_{\text{CR}}$ | $\text{GAP}_{\text{BB}}$ |
| 3 | 0.1 | 0.5 | 7.88 | 0.03 | 10 | 0 | 56.42 | 0.07 | | 1 | 0 | 3405.19 | 93.64 | 36.33 |
| | 0.5 | 0.5 | 556.76 | 0.00 | 2 | 2 | 2359.87 | 0.07 | 0.34 | 0 | 10 | † | | |
| 5 | 0.1 | 0.5 | 1.41 | 0.14 | 10 | 0 | 3.86 | 0.24 | | 4 | 0 | 1329.74 | 89.57 | 33.05 |
| | 0.1 | 0.7 | 2.87 | 0.32 | 10 | 0 | 5.53 | 0.51 | | 7 | 0 | 1709.50 | 89.59 | 35.02 |
| | 0.3 | 1 | 528.39 | 0.00 | 8 | 0 | 2608.07 | 0.09 | 0.33 | 0 | 1 | † | | 36.05 |
| | 0.5 | 0.5 | 168.13 | 0.00 | 10 | 0 | 926.94 | 0.09 | | 0 | 0 | † | | 35.29 |
| 7 | 0.1 | 0.7 | 1.05 | 0.48 | 10 | 0 | 3.83 | 0.76 | | 9 | 0 | 2477.41 | 85.73 | 33.81 |
| | 0.3 | 1 | 175.48 | 0.00 | 10 | 0 | 457.61 | 0.18 | | 0 | 0 | † | | 36.12 |

#O = number of opt. sol.; #N = number of instances without sol. within time limit
CPU='†' if time limit reached (3600 s).
$\text{DEV}_{\text{CR}}$=' ' if non opt or no sol.; $\text{GAP}_{\text{BB}}$=' ' for opt. sol. or no sol.

Table 5: Performance of FFT, FFF and TIV I for $n = 100$ jobs

The high performance of our formulations can also be explained by the reduced number of variables they allow for, compared with TIV I as evidenced in Figure 2. We also reported in parenthesis the number of optimal solutions reached by TIV I for each combination of parameters (indicated solely by the number of machines on the $x-$axis). Although we can not state a clear relationship between the number of machines and the number of variables, neither between the number of variables of TIV I and its optimality rate, it strikingly appears that FFT dramatically reduces the number of variables, hence the search space.

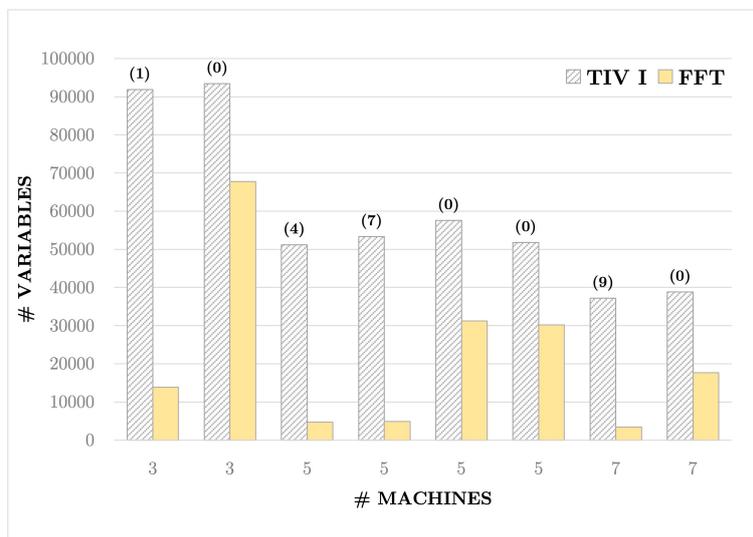

Figure 2: Average number of variables of FFT and TIV I for $n = 100$ jobs

For large-sized instances with more than 100 jobs, due to the inability of TIV I to find integer solutions in the time limit of 3600 s, we only considered FFT and FFT-Warmed, a variant of FFT with initial solutions as provided by heuristics HS1, HS2 of Elidrissi et al. (2018b). Table 6 displays the results for $n = 150, 200$. Easy instances (parameters in black font, 120 instances out of 240) were solved to optimality by FFT in 7.22 s and 15.48 s for $n = 150$ and $n = 200$ jobs respectively and FFT-Warmed allowed for a gain in execution time of about 40% on average. But the benefit of using initial solutions, namely FFT-Warmed, was even clearer for some hard instances (bold blue font) for which FFT even failed to find integer solutions for some draws within the time limit. For these hard parameters combinations, FFT-Warmed produced solutions to all the draws with a very low average $\text{GAP}_{\text{BB}}$ (less than 1.33% at maximum). It should be noted that the hardest instances with $(\alpha, \rho) = (0.5, 0.5)$ for $n = 150$ could be solved to optimality by FFT-Warmed in less than 3 hours on average (2.83 hours). Again, we can observe from Table 6 that the continuous relaxation of our Flow-Flow model FFT generated lower bounds of high quality, with $\text{DEV}_{\text{CR}}$ all less than 0.28% on average (obviously $\text{DEV}_{\text{CR}}$ are the same for FFT-Warmed).

| | | | | **FFT** | | | | | **FFT-Warmed** | | |
|---|---|---|---|---|---|---|---|---|---|---|---|
| $n$ | $m$ | $\alpha$ | $\rho$ | #O | #N | CPU | $\text{DEV}_{\text{CR}}$ | $\text{GAP}_{\text{BB}}$ | #O | CPU | $\text{GAP}_{\text{BB}}$ |
| 150 | 3 | 0.1 | 0.5 | 10 | 0 | 17.83 | 0.01 | | 10 | 8.38 | |
| | | **0.5** | **0.5** | 3 | 7 | 2875.46 | 0.00 | | 5 | 2424.28 | 0.34 |
| | 5 | 0.1 | 0.5 | 10 | 0 | 2.82 | 0.08 | | 10 | 2.45 | |
| | | 0.1 | 0.7 | 10 | 0 | 5.89 | 0.23 | | 10 | 4.62 | |
| | | 0.3 | 1 | 10 | 0 | 2072.49 | 0.00 | | 10 | 1673.31 | |
| | | 0.5 | 0.5 | 10 | 0 | 840.60 | 0.00 | | 10 | 595.98 | |
| | 7 | 0.1 | 0.7 | 10 | 0 | 2.34 | 0.28 | | 10 | 2.23 | |
| | | 0.3 | 1 | 10 | 0 | 498.60 | 0.00 | | 10 | 388.20 | |
| 200 | 3 | 0.1 | 0.5 | 10 | 0 | 39.85 | 0.01 | | 10 | 14.16 | |
| | | **0.5** | **0.5** | 0 | 10 | † | | | 0 | † | 0.24 |
| | 5 | 0.1 | 0.5 | 10 | 0 | 3.84 | 0.07 | | 10 | 2.91 | |
| | | 0.1 | 0.7 | 10 | 0 | 12.58 | 0.16 | | 10 | 7.71 | |
| | | **0.3** | **1** | 2 | 6 | 2866.58 | 0.00 | 0.17 | 0 | † | 1.33 |
| | | 0.5 | 0.5 | 10 | 0 | 2215.44 | 0.00 | | 8 | 2710.35 | 0.64 |
| | 7 | 0.1 | 0.7 | 10 | 0 | 5.64 | 0.26 | | 10 | 5.57 | |
| | | 0.3 | 1 | 10 | 0 | 1719.48 | 0.00 | | 10 | 1008.92 | |

#O = number of opt. sol.; #N = number of instances without sol. within time limit

CPU='†' if time limit reached (3600 s).

$\text{DEV}_{\text{CR}}$=' ' if non opt or no sol.; $\text{GAP}_{\text{BB}}$=' ' for opt. sol. or no sol.

Table 6: Performance of FFT and FFT-Warmed for $n = 150, 200$ jobs

## 4 Conclusion

In this paper, we have proposed an arc-flow formulation to solve the scheduling problem of a set of jobs that must be processed non preemptively on identical parallel machines and requiring prior setup operations on a common server, with the aim of minimising the makespan. Our model relies on an arc-flow formulation for both machines and server with shared variables related to the start time of the jobs and uses an efficient procedure for handling identical jobs. The resulting Flow-Flow formulation (FFF) therefore works with a fairly limited number of variables. Furthermore, the continuous relaxation of FFF produces high quality lower bounds that expedite the convergence towards optimal solutions. By contrast, the best formulation in past research, TIV I, often generates lower bounds below a trivial hence poor lower bound to the problem. Computational experiments showed that FFF was very effective in solving problems with up to 100 jobs and 7 machines, especially its tuned version FFT whereas TIV I failed to find optimal solutions in 74% of the cases. Extended experiments with $n = 150, 200$ jobs again highlighted the good performance of our model, with 83% of optimal solutions and integer solutions in the time limit with extremely low gaps to best bound as provided by CPLEX. Our experimental framework evidenced hard instances to the scheduling problem at hand, which suggests that additional effort could be spent to enhance our solution method, with column generation as a promising avenue.

Besides, our Flow-Flow formulation could be adapted to several variants of the problem that appear in real settings such as configurations with two identical servers or a single server performing pre and post operations.